\documentclass[10pt,a4paper,twoside,leqno]{article}
\usepackage{amsfonts,amssymb}
\usepackage{eufrak}
\usepackage{amscd}

\font\Bigtit=cmr10 scaled \magstep 4
\font\ebf=cmbx8
\font\erm=cmr8

\usepackage{graphicx}

\begin{document}

\thispagestyle{empty}

\begin{flushright}
PL ISSN 0459-6854
\end{flushright}
\vspace{0.5cm}
\centerline{\Bigtit B U L L E T I N}
\vspace{0.5cm}
\centerline{DE \ \  LA \ \  SOCI\'ET\'E \ \  DES \ \  SCIENCES \ \ ET \ \ DES \
\ \ LETTRES \ \ DE \ \ \L \'OD\'Z}
\vspace{0.3cm}
\noindent 2003\hfill Vol. LIII
\vspace{0.3cm}
\hrule
\vspace{5pt}
\noindent Recherches sur les d\'eformations \hfill Vol. XLII
\vspace{5pt}
\hrule
\vspace{0.3cm}
\noindent pp.~27--37

\vspace{0.7cm}

\noindent {\it Andrzej K. Kwa\'sniewski}

\vspace{0.5cm}

\noindent {\bf ON FIBONOMIAL AND OTHER TRIANGLES  VERSUS  DUALITY TRIADS}

\vspace{0.5cm}

\noindent {\ebf Summary}

{\small
The duality  triads were defined in the preceding paper. Notation,
enumeration of formulas  and references is therefore to be continued
hereby. In this paper Fibonomial   triangles [3--8] and further Pascal-like
triangles including the $q$-Gaussian one are given explicit interpretation as
discrete time dynamical systems as it is in the case with duality triads.
Large classes of duality triads in various areas of  mathematics  are
identified including those triads that appear in finite operator calculus
and its extensions [9--17].}

\vspace{0.5cm}

\renewcommand{\thesubsection}{\arabic{subsection}.}

\subsection{New examples of dual triads  and related\\ general information}

The duality  triads were defined in the preceding paper [22].

\subsubsection{$q$-Gaussian dual triad}

(see: pp. 68-74, Proposition 12.1, recurrence 12.4,  Propositions 12.3
and 12.4 in [32]). Let us recall that $q$-Gaussian [32, 33, 34,16]
coefficient is given by
$$
\left( \begin{array}{c} n\\k \end{array}
\right)=\frac{n_{q}!}{k_{q}!(n-k)_{q}!}=
\left( \begin{array}{c} n\\n-k\end{array} \right)_{q}.
 $$

\noindent  For  a positive integer  $q >1; q = p^{m}$  where  $p$  is
prime $ \left( \begin{array}{c} n\\k \end{array}
\right)_{q}$  is interpreted  [32, 35] as  the number of $k$-dimensional
subspaces in n-th dimensional space over Galois field   $GF(q)$.  Also  $q$
real  and  $-1<q<+1$ are exploited in vast literature  (see: (2.5) in [33]
and references therein).  Here for
$$n>0, \quad n_{q}\equiv \frac{1-q^{n}}{1-q}$$

\noindent     and while  using the
following notation for $q$-factorial (see also  [34, 17]):
$$
n_{q}!=n_{q}(n-1)_{q}!;\quad 1_{q}!=0_{q}!=1,\quad
n_{q}^{\underline{k}}=n_{q}(n-1)_{q}\ldots (n-k+1)_{q}, $$

\noindent we write
$$
\left( \begin{array}{c} n\\k \end{array}
\right)_{q} \equiv \frac{n_{q}^{\underline{k}}}{k_{q}!}.
$$

\noindent The recurrence relation   [16, 32, 33] for connection constants
$$
C_{n,k}=\left( \begin{array}{c} n\\k \end{array}
\right)_{q}
$$

\noindent reads as follows:
\begin{eqnarray}
&&\left( \begin{array}{c} n+1\\k \end{array}
\right)_{q} =q^{k}\left( \begin{array}{c} n\\k \end{array}
\right)_{q} + \left( \begin{array}{c} n\\k-1 \end{array}
\right)_{q}, \nonumber \\
&&\label{eq*7}\\
&& \left( \begin{array}{c} n\\0 \end{array}
\right)_{q} =1, \quad n \geq 0,\quad k\geq 1 \nonumber
\end{eqnarray}

\noindent and might be considered as the result of the $q$-Leibniz rule in
q-umbral calculus\linebreak (p. 176 in [16]) of Eulerian formal power series
[12, 32, 35]. Another  equivalent form  is the following
$$
\left( \begin{array}{c} n+1\\k \end{array}
\right)_{q} =\left( \begin{array}{c} n\\k \end{array}
\right)_{q}  +q^{n-k}\left( \begin{array}{c} n\\k -1 \end{array}
\right)_{q} ,
$$
$$
\left( \begin{array}{c} n\\0 \end{array}
\right)_{q} =1, \quad n \geq 0,\quad k \geq 1.
$$

\noindent (The Eulerian formal power series  as well as many others form
algebras isomorphic to the reduced incidence algebras  [12, 32, 35]).  The
above recurrences (\ref{eq*7}) may be represented by the  $q$-Gaussian
triangles. The corresponding  $q=2,\ q=3,\  q = 5$  cases are supplied below:

\noindent (I) $q=2$
{\setlength\arraycolsep{2pt}
$$
\begin{array}{ccccccccccccccc}
&&&&&&&{1}&&&&&&&\\
&&&&&&{1}&&{1}&&&&&&\\
&&&&&{1}&&3&&{1}&&&&&\\
&&&&{1}&&7&&7&&{1}&&&&\\
&&&{1}&&{15}&&35&&15&&{1}&&&\\
&&{1}&&31&&{155}&&155&&31&&{1}&&\\
&{1}&&63&&{651}&&1395&&651&&{63}&&{1}&\\
\bullet&&\bullet&&\bullet&&\bullet&&\bullet&&\bullet&&\bullet&&\bullet
\end{array} $$}

\noindent (II) $q=3$
{\setlength\arraycolsep{2pt}
$$
\begin{array}{ccccccccccccccc}
&&&&&&&{1}&&&&&&&\\
&&&&&&{1}&&{1}&&&&&&\\
&&&&&{1}&&4&&{1}&&&&&\\
&&&&{1}&&13&&13&&{1}&&&&\\
&&&{1}&&{40}&&130&&40&&{1}&&&\\
&&{1}&&121&&{1210}&&1210&&121&&{1}&&\\
&{1}&&364&&{11011}&&3388&&11011&&{364}&&{1}&\\
\bullet&&\bullet&&\bullet&&\bullet&&\bullet&&\bullet&&\bullet&&\bullet
\end{array} $$}

\noindent (III) $q=5$
{\setlength\arraycolsep{2pt}
$$
\begin{array}{ccccccccccccccc}
&&&&&&&{1}&&&&&&&\\
&&&&&&{1}&&{1}&&&&&&\\
&&&&&{1}&&6&&{1}&&&&&\\
&&&&{1}&&31&&31&&{1}&&&&\\
&&&{1}&&{156}&&806&&156&&{1}&&&\\
&&{1}&&781&&{20306}&&20306&&781&&{1}&&\\
&{1}&&3906&&{508431}&&16401&&508431&&{3906}&&{1}&\\
\bullet&&\bullet&&\bullet&&\bullet&&\bullet&&\bullet&&\bullet&&\bullet
\end{array} $$}

\noindent The   recurrence dual to (\ref{eq*7}) is then  given by
\begin{eqnarray}
&&n\Phi_{n}(x)=q^{n}\Phi_{n}(x)+\Phi_{n+1}(x),\nonumber\\
&&\label{eq**7}\\
&&\Phi_{0}(x)=1,\quad \Phi_{-1}(x)=0,\quad; \quad n \geq 0.\nonumber
\end{eqnarray}

\noindent in accordance with a well known fact [32]  that
\begin{equation}\label{eq***7}
x^{n}=\sum\limits _{k \geq 0}\left( \begin{array}{c} n\\k
\end{array} \right)_{q}\Phi_{k}(x), \end{equation}
$$
\Phi_{k}(x)=\prod\limits_{s=0}^{k-1} \left( x-q^{s}  \right).
$$

\noindent $\Phi_{k}(x)=\prod\limits_{s=0}^{k-1} \left( x-q^{s}  \right)$
are named to be  $q$-Gaussian polynomials  and  note that these are monic
persistent root polynomials [26]  therefore the recurrence  (\ref{eq*7}) is
of the form
$$
[r]=\left\{\left( q^{k-1} \right)  \right\}_{k \geq 1},
$$
\begin{eqnarray}
&&L_{n+1,k}=L_{n,k-1}+r_{k+1}L_{n,k},\nonumber \\
&&\label{eq-qG}\\
&&L_{0,0}=1,\quad L_{0,-1}=0,\quad n,k \geq 0, \nonumber
\end{eqnarray}

\noindent in the notation of  [26]  i.e.
$$
\left( \begin{array}{c} n\\k
\end{array} \right)_{q}  =c_{n,k}=L_{n,k}.
$$

\subsubsection{Catalan triad  and triangle}

Let  $C_{n,k}$  denotes the number of pairs of nonintersecting paths of
length $n$  and distance  $k$  as defined in [36].  Then one may prove
[36] that $C_{n,k}$ (where now $n \geq k >0$)  satisfy the  recurrence
\begin{eqnarray}
&&C_{n+1,k}=C_{n,k-1}+2C_{n,k}+C_{n,k+1},\nonumber\\
&&\label{eq*C}\\
&&C_{1,1}=1,\quad C_{n,0}=0=C_{n,n+k}=0 \quad {\rm for}\quad n \geq k\geq
0.\nonumber \end{eqnarray}

\noindent The transition matrix is then of the form
$$
X =\left( \begin{array}{ccccccc} 2&1&0&0&0&0&\ldots\\
1&2&1&0&0&0&\ldots\\
0&1&2&1&0&0&\ldots\\
0&0&1&2&1&0&\ldots\\
0&0&0&1&2&1&\ldots  \end{array} \right)
$$

\noindent and we obtain the Catalan triangle  calculating its subsequent
rows  $C_{n} = C_{0} X^{n};  n \geq 0$ in accordance with the fact [36]
that
$$
C_{n,k}=\left( \begin{array}{c} 2n\\n-k
\end{array} \right)\frac{k}{n}:
$$

\noindent (note:
$$
C_{n,1}=\left( \begin{array}{c} 2n\\n-1
\end{array} \right) \frac{1}{n}=\left( \begin{array}{c} 2n\\n
\end{array} \right) \frac{1}{n+1}
$$

\noindent are Catalan numbers).

The Catalan  Triangle
$$
\begin{array}{ccccccccccccc}
&&&&&&{1}&&&&&&\\
&&&&&{0}&&{1}&&&&&\\
&&&&{0}&&2&&{1}&&&&\\
&&&{0}&&5&&4&&{1}&&&\\
&&{0}&&{14}&&14&&6&&{1}&&\\
&{0}&&42&&{48}&&27&&8&&{1}&\\
\bullet&&\bullet&&\bullet&&\bullet&&\bullet&&\bullet&&\bullet \end{array}
$$

\noindent with dual  recurrence for  the corresponding polynomials being of the form
\begin{eqnarray}
&&xC_{k}(x)=C_{k+1}(x)+2C_{k}(x)+C_{k-1}(x),\quad k \geq , 0\nonumber \\
&&\label{eq**C}\\
&&C_{0}(x)=1, \quad C_{-1}(x)=0.\nonumber
\end{eqnarray}

\noindent To complete the triad we note according to general rule  that
\begin{equation}\label{eq***C}
x^{n}=\sum\limits _{1\leq k\leq n}  \left( \begin{array}{c} 2n\\n-k
\end{array} \right) \frac{k}{n}C_{k}(x),\quad n>0.
\end{equation}

\noindent We call the polynomial sequence $\{C_{k}(x)\}_{k \geq 0}$  the
Catalan polynomial sequence.

\subsubsection{Remarks and Information IV}

The $q$-Gaussian polynomials
$$
\Phi_{k}(x)=\prod\limits_{s=0}^{k-1}\left( x-q^{s} \right)
$$

\noindent  from Example  7 (along with some
other examples)  constitute a  persistent root polynomial sequence and
therefore as such they do satisfy the recurrence equation  [23]   (see:
Remark and Information  I)
\begin{eqnarray}
&&L_{n+1,k}=L_{n,k-1}+r_{k+1}L_{n,k}, \nonumber \\
&&\label{eq*L}\\
&&L_{0,0}=1,\quad L_{0,-1}=0 \quad n,k\geq 0, \nonumber
\end{eqnarray}

\noindent where  $[r] \equiv \{ r_{0}, r_{1}, r_{2},\ldots ,
r_{k},\ldots \}$ is the root sequence determining $\{\Phi_{n} \}_{n \geq
0}$ and the connection constants $(L_{n,k})$ -- called in [23] the
generalized Lah numbers -- are determined by duality triad  (\ref{eq***C})
equation
\begin{equation}\label{eq***Fi}
x^{n}=\sum\limits_{0\leq k\leq n}L_{n,k}\Phi _{k}(x),\quad n \geq 0 ,
\end{equation}

\noindent thus completing another example of duality triad. Naturally this
is the case with any
persistent root   polynomial sequence $\{\Phi _{n}  \}_{n \geq 0}$. On
this important occasion note that because of Farward Theorem  (see: p. 21 in
[29] and [26]) persistent root polynomial sequences do not form OPS
sequences  (OPS = Orthonormal Polynomial System -- [29]) because OPS
sequences should satisfy three term recurrences. An example of OPS  at hand
is provided by Catalan polynomials satisfying (\ref{eq**C}). The other
examples are classical Hermite or Laguerre polynomials  or Tchebychev OPS
(in this last case -- $X$ represents  the symmetric Jacobi matrix see:
p.~138 in [29]).

\vspace{3mm}

\noindent {\it Concluding}:  persistent root polynomials sequences $\{ \Phi
_{n} \}_{n \geq 0}$ and (disjointedly)  OPS sequences -- orthogonal  with
respect to quasi-definite moment functional -- give rise to corresponding
duality triads. Some of  OPS sequences like  Hermite, Laguerre,
$q$-Hermite [16, 37],  $q$-Laguerre (see: (33) in [38]  and  [17])    etc.
are generalized Appell polynomials.

\vspace{3mm}

\noindent {\it Note 1}: not all Appell polynomials give rise to duality
triads (see: Remarks and Information I).

\vspace{3mm}

\noindent {\it Note 2}: all  generalized Appell polynomials that are
at the same time persistent root polynomials  are discovered by the authors
of  [26] and these polynomial sequences -- naturally --  give rise  to
duality triads.

\subsection{Does Fibonomial  dual triad  exist? }

Let us start investigating this problem with reference to   $q$-Gaussian
dual triad  example. The corresponding recurrence relation   connection constants
$c_{n,k}$   might be considered as the result of the $q$-Leibniz rule in
$q$-umbral calculus (p.~176 in [16]) of Eulerian formal power series
[12, 32, 35], where that series  as well as many others form
algebras isomorphic to the standard reduced incidence algebras  [12, 32, 35].
The connection constants   $c_{n,k}=\left( \begin{array}{c} n\\k
\end{array} \right)_{q}$ from relation (\ref{eq***Fi}) of the triad under
consideration are interpreted [12, 32, 35] as  the number of $k$-dimensional
subspaces in $n$-th dimensional space $V(n,q)$ over the Galois field $GF(q)$.
Restating this in incidence algebras language [35]  we conclude that in the
lattice $L(n,q)$ of all subspaces of $V(n,q)$:

\begin{center}
$n_{q}!=$ {\it the number of maximal chains in an of  interval  of length
$n$,}

$$\left( \begin{array}{c} n\\k
\end{array} \right)_{q}=\left[ \begin{array}{c} \alpha \\ \beta \gamma
\end{array} \right]:=$$
{\it the number of distinct elements  $z$  in a
segment $[x,y]$ of type $\alpha$  such  that}

$[x,z]$ is of type $\beta$ while $[z,y]$ is of type $\gamma$,
\end{center}

\noindent where  $x,y,z  \in L(n,q)$  and the type of  $[x,y]$ is equal to
${\rm dim}[x/y]$. (Here  $[y/x]$ denotes the subtraction of linear
subspaces [35]).

In straightforward analogy consider now {\it  Fibonomial  coefficients}
[3--8], [39--41]   $c_{n,k}$, given by
 $$
 c_{n,k}= \left( \begin{array}{c} n\\k
\end{array} \right)_{F}= \frac{F_{n}!}{F_{k}!F_{n-k}!}=\left(
\begin{array}{c} n\\n-k \end{array} \right)_{E}, $$

\noindent where  [34, 43]:
$$
n_{F}\equiv F_{n} \neq 0,\quad
n_{F}!\equiv n_{F}(n-1)_{F}(n-2)_{F}(n-3)_{F}\ldots 2_{F}1_{F};\quad
0_{F}!=1; $$
$$
n_{F}^{\underline{K}}=n_{F}(n-1)_{F}\ldots (n-k+1)_{F}; \quad \left( \begin{array}{c} n\\k
\end{array} \right)_{F} \equiv \frac{n_{F}^{\underline{k}}}{k_{F}!}.
$$

\noindent A natural question then arises:

\vspace{3mm}

{\it Does there exist a partialy ordered set   $P$  with its  incidence
algebra $R (P)$  and  such incidence coefficients  $\left[ \begin{array}{c}
\alpha \\ \beta \gamma\end{array}\right]$ that coincide with $F$-binomial
ones i.e. $\left[ \begin{array}{c}
\alpha \\ \beta \gamma\end{array}\right]=\left( \begin{array}{c}
n \\ k\end{array}\right)_{F}$}?

\vspace{3mm}

For the moment this is an open problem for us. However independently of
that  one may naturally consider  -- in analogy to Eulerian formal series
algebra -- also a  ``Fibonomial'' formal series algebra with a Fibonomial
convolution as the associative product  [43, 26, 48], i.e.
$$
c_{n}=\sum\limits_{k \geq 0}^{n} \left( \begin{array}{c}
n \\ k\end{array}\right)_{F} a_{k}b_{n-k},\quad {\rm where}\quad
F(z)G(z)=H(z)\equiv \sum\limits_{n \geq 0}\frac{c_{n}}{n!}z^{n}, $$
$$
G(z)\equiv \sum\limits_{n \geq 0}\frac{b_{n}}{n!}z^{n}\quad {\rm and}\quad
F(z)\equiv \sum\limits_{n \geq 0} \frac{a_{n}}{n!}z^{n}. $$

The corresponding Fibonomial finite operator calculus [43] becomes then an
example  of the general theory of the ``calculus of sequences'' [34]  in
its operator form [13--15], [49]. As for the recurrence (see: p.~27,formule
(58) in [41]) for  {\it Fibonomial  coefficients}  $c_{n,k}$, it  is well
known [39--41] ([42] is also relevant)  and some properties  of  Fibonomial
numbers  $\left(\begin{array}{c} n\\k \end{array}  \right)_{F}$  i.e. of
Fibonomial triangle were among others considered in [3--8], [40], [44--48].
This recurrence reads:
\begin{eqnarray}
&&\left(\begin{array}{c} n+1\\k \end{array}  \right)_{F} =F_{k-1}
\left(\begin{array}{c} n\\k \end{array}
\right)_{F}+F_{n-k+2}\left(\begin{array}{c} n\\k-1 \end{array}  \right)_{F} ,
\nonumber \\
&&\label{eq*F}\\
&&\left(\begin{array}{c} n\\0 \end{array}  \right)_{F}=1,\quad
\left(\begin{array}{c}0\\k \end{array}  \right)_{F} =0
\quad{\rm for}\quad k>0 \nonumber
\end{eqnarray}

\noindent or, equivalently,
\begin{eqnarray*}
&&\left(\begin{array}{c} n+1\\k \end{array}  \right)_{F} =F_{k+1}
\left(\begin{array}{c} n\\k \end{array}
\right)_{F}+F_{n-k}\left(\begin{array}{c} n\\k-1 \end{array}  \right)_{F},
\\
&&\\
&&\left(\begin{array}{c} n\\0 \end{array}  \right)_{F}=1,\quad
\left(\begin{array}{c}0\\k \end{array}  \right)_{F} =0
\quad{\rm for}\quad k>0 .
\end{eqnarray*}

\noindent Note that the coefficients of  the recurrence  (\ref{eq*F})  for
$c_{n,k}$ depend on both  $k$  and $n$ as it was not the case in all
previous examples. We face now  a new situation out of the scope of the
former definition of dual triad, where it was assumed that in
\begin{eqnarray}
&&c_{n+1,k}=i_{k-1}c_{n,k-1}+q_{k}c_{n,k}+d_{k+1}c_{n,k+1}\nonumber\\
&& \label{eq*0}\\
&&c_{0,0}=1;\quad c_{0,k}=0\quad {\rm for}\quad k>0. \nonumber
\end{eqnarray}

\noindent the numbers $i_{k},\ q_{k},\ d_{k};\ k \geq 0$ -- are
{\it independent} of ``discrete time'' parameter  $n$.

This assumption
was decisive in deriving  the third completing member of
the triad. Nevertheless -- it seems -- we still may maintain a ``discrete
time'' dependent interpretation. Namely $c_{n,k} =$ {\it number of ways to
reach level $k$ in $n$ steps starting from the level} 0 and still
(\ref{eq*F}) may be described  as such a discrete time system in which
corresponding rows (``states'') of the $C = (c_{n,k})$  matrix   (hence:
rows of Fibonomial triangle) are reached via  consecutive application of
different, appropriately adjusted, unique one step  transition matrix $F$,
which however is not that of dual recurrence -- simply, because dual
recurrence does not exist. Nevertheless,  generalizing [1]  we still may
maintained the standpoint
the  $kl$-th entry of  $F(n)  =$ {\it number of ways of going from the level
$k \to to \to l$   in $n$ steps}   where $(c_{n,k} \neq 0 \ {\rm for}\
k\leq n )\ F$ is the solution of recurrent for rows of  $C$ equation  $C F
= EC$ from which it follows immediately that $E^{n} C = C F^{n},\  n \geq 0$,
where
$$
X =\left( \begin{array}{ccccc} 0&1&0&0&\ldots\\
0&0&1&0&\ldots\\
0&0&0&1&\ldots\\
0&0&0&0&\ldots\\
\ldots&\ldots&\ldots&\ldots&\ldots  \end{array} \right) =\left( \delta
_{i+1,k} \right)_{i,k \geq 0}. $$

\noindent Such a matrix  $F$  exists  -- in our case -- giving rise in a
standard way to the Fibonomial triangle
{\setlength\arraycolsep{2pt}
$$
\begin{array}{ccccccccccccccc}
&&&&&&&{\bf 1}&&&&&&&\\
&&&&&&{\bf 1}&&{\bf 1}&&&&&&\\
&&&&&{\bf 1}&&1&&{\bf 1}&&&&&\\
&&&&{\bf 1}&&2&&2&&{\bf 1}&&&&\\
&&&{\bf 1}&&{\bf 3}&&{\bf 6}&&3&&{\bf 1}&&&\\
&&{\bf 1}&&5&&{\bf 15}&&15&&5&&{\bf 1}&&\\
&{\bf 1}&&8&&{40}&&60&&40&&{8}&&{\bf 1}&\\
\bullet&&\bullet&&\bullet&&\bullet&&\bullet&&\bullet&&\bullet&&\bullet \end{array}
$$}

\noindent which, as said, might be obtained  when considering subsequent
powers of  this triangle generating matrix $F$
$$
X =\left( \begin{array}{ccccccc} 1&1&0&0&0&0&\ldots\\
0&0&1&0&0&0&\ldots\\
0&1&1&1&0&0&\ldots\\
0&0&2&1&1&0&\ldots\\
0&-2&0&6&2&1&\ldots  \end{array} \right)
$$

\noindent and then   calculating  $C_{n} = C_{0} F^{n};\   n \geq 0$,   where  $C_{0}
= (1 , 0 , 0 , 0 , . . . )$. The subsequent rows of  $F$  are:
\begin{eqnarray*}
&&0,\ 2,\ -10,\ 0,\ 15,\  3,\   1,\ 0,\     0,\  0,\ 0,\  0,\   0, \ldots\\
&&0,\ 36,\ 16,\ -80,\ -100,\ 40,\ 5,\    1,\   0,\  0,\ 0,\  0,\ldots\\
&&0,\ \ldots,\ \ldots ,\ \ldots ,\  \ldots ,\  \ldots ,\ \ldots ,\  8,\
1,\ 0,\  0,\ 0,\ldots\\
&&\ldots\ldots\ldots\ldots\ldots\ldots\ldots\ldots\ldots\ldots\ldots\ldots\ldots\ldots\ldots\ldots
\end{eqnarray*}

One here just adjusts  $F$ matrix elements step by step while
calculating $ C_{n} = C_{0} F^{n}$ The formula for matrix elements of  $F$
is beyond our abilities and needs. The formula for the $k$-th component
$[C_{n}]_{k}$ of $C_{n}$  is tautologically obvious:
$$
\left[C_{n} \right]_{k} =\left( \begin{array}{c} n\\k \end{array} \right)
_{F}. $$

\noindent Eventually expected recurrence in the form  of an eigenvector and
eigenvalue equation (6)   $x\Phi = F\Phi$  by no means would be dual
recurrence  to the recurrence (\ref{eq*F}). The resulting sequence $\Phi$
is not a polynomial sequence $({\rm deg}\ p_{n} = n)$. We may see it
immediately from the start:
\begin{eqnarray*}
&&\Phi_{0}(x) = 1,\quad \Phi_{1}(x) = x-1,\quad  \Phi_{2}(x) = 0,\quad
\Phi_{3}(x) =1-x,\\
&& \Phi_{4}(x) = -(x-1)^{2},\quad \Phi_{5}(x) =
(2-x)(x-1)^{2} +8(x-1).\end{eqnarray*}

\noindent Moreover  any polynomial sequence solution  would be dual
recurrence:
\begin{eqnarray}
&&x\Phi_{k}(x)=d_{k}(n)\Phi_{k-1}(x)+q_{k}(n)\Phi_{k}(x)+i_{k}(n)\Phi_{k+1}(
x); \nonumber \\
&&\label{eq**0prime}\\
&&\Phi_{0}(x)=1,\quad \Phi_{-1}(x)=0;\quad k \geq 0,\nonumber \end{eqnarray}

\noindent which does not exist for numbers $i_{k},\ q_{k},\ d_{k};\ k \geq 0$ --
{\it dependent} of ``discrete time'' parameter  $n$ as it is the case with the
Fibonomial coefficients recurrence.

\vspace{3mm}

\noindent {\it Conclusion}.
Neither  dual recurrences nor triad polynomials  exist in the sense of numbers for
 $i_{k},\ q_{k},\ d_{k};\ k \geq 0$  being {\it dependent} of
``discrete time'' parameter $n$. The Fibonomial case under consideration is
one of such examples as
\begin{equation}\label{eq13}
c_{n,k}=\left[ \begin{array}{c} n\\k \end{array}
\right], \quad c_{n,k}=\left\langle \begin{array}{c} n\\k \end{array}
\right\rangle, \quad c_{n,k}=\left( \begin{array}{c} n\\k \end{array}
\right)_{F}.
\end{equation}

\subsection{Families of  triad polynomials }

Existence of dual triad is the quite specific property of corresponding
objects as there are many  polynomial sequences of distinguished importance
which are  triad polynomial sequences as well as many polynomial sequences
of distinguished importance which are not triad polynomial sequences:  for
example Abel, Euler  or exponential polynomials.

Here we summarize the information supplied above as to {\it when} a family
of polynomial sequences   $\{ \Phi _{n} \}_{n \geq 0}({\rm deg} \Phi
_{n}(x) = n)$ is a family of triad polynomials.
\begin{enumerate}
\item First of all let us recall that     (\ref{eq*0}) and
(\ref{eq**0prime}) imply (\ref{eq13})  but not vice versa as (\ref{eq13})
does not imply (\ref{eq*0}) (Remarks and Information I).

\item Secondly for numbers $i_{k},\ q_{k},\ d_{k};\ k \geq 0$ --
{\it dependent} of  ``discrete time'' parameter  $n$  in the recurrence
(\ref{eq*0}) triads do not exist. Famous examples are:
$$ c_{n,k}=\left[ \begin{array}{c} n\\k \end{array}
\right], \quad c_{n,k}=\left\langle \begin{array}{c} n\\k \end{array}
\right\rangle, \quad c_{n,k}=\left( \begin{array}{c} n\\k \end{array}
\right)_{F}. $$

\noindent The polynomial sequence solution would be dual recurrence
(\ref{eq**0prime}) which does not exists for numbers $i_{k},\ q_{k},\ d_{k};\ k \geq
0$ -- dependent of ``discrete time'' parameter  $n$.

\item Any persistent roots polynomial sequence is a sequence of triad
polynomials (All  generalized Appell polynomials that are at the same time
persistent root polynomials are given in  [26]).

\item Any OPS  -- with respect to quasi-definite moment functional -- is a
sequence of triad polynomials.

\item All OPS being at the same time  generalized   Appell polynomial
sequences have been determined by Chihara in [50] (see also [29]  p.~167
and [51] and [52], [9, 12, 29], [53]).
\end{enumerate}






\vspace{0.5cm}

\noindent {\erm Institute of Computer Science}

\noindent {\erm Bia\l ystok University }

\noindent{\erm  Sosnowa 64, PL-15-887 Bia\l ystok}

\noindent {\erm Poland}

\vspace{0.5cm}

\noindent Presented by Julian \L awrynowicz at the Session of the
Mathematical-Physical Commission of the \L \'od\'z Society of Sciences and
Arts on December 17, 2003

\vspace{0.5cm}
\noindent {\bf POR\'OWNANIE TR\'OJK\c{A}T\'OW FIBONOMIALNYCH \\I POKREWNYCH
Z TR\'OJKAMI DUALNYMI}

\vspace{0.2cm}

\noindent {\small S t r e s z c z e n i e}

{\small W pracy  poprzedzaj\c{a}cej wprowadzono poj\c{e}cie triad dualnych,
o kt\'orych pe\l na informacja mo\.ze by\'c przedstawiona w postaci
``tr\'ojk\c{a}ta'' na podobie\'nstwo tr\'ojk\c{a}ta Pascala.
W niniejszej pracy b\c{e}d\c{a}cej kontynuacj\c{a} poprzedniej rozwa\.za
si\c{e} tr\'ojk\c{a}t fibonomialny\linebreak[4] [3--8], $q$-Gaussowski i inne
interpretowane jako uk\l ady dynamiczne z czasem dyskretnym.
W~pracy zidentyfikowano dwie  obszerne klasy  triad  reprezentowane
poprzez  ci\c{a}gi wielomian\'ow quasi-ortogonalnych (w tym ortogonalnych)
[29] oraz te reprezentowane przez  ci\c{a}gi wielomian\'ow znanych pod
nazw\c{a} ``persistent roots polynomials''  [26].
Po\'sr\'od tych przyk\l ad\'ow znajdujemy interesuj\c{a}ce przypadki z
klasycznego i rozsze\-rzonego sko\'nczonego rachunku operatorowego Roty
[12].}



\end{document}